\documentclass[11pt]{article}

\usepackage{amscd,amsmath, amssymb, fancyhdr,
  epsfig,url,MnSymbol}
\usepackage[backref=page]{hyperref}

\usepackage{color}
\usepackage[matrix,arrow]{xy}
\usepackage{pb-diagram}
\numberwithin{equation}{section}

\renewcommand*{\backrefalt}[4]{%
	\ifcase #1 (Not cited.)%
	\or        (Cited on page~#2.)%
	\else      (Cited on pages~#2.)%
	\fi}

\hypersetup{
	colorlinks   = true,
	citecolor    = black,
	linkcolor    =black,
	urlcolor     =black	
}

\newcommand{\version}{July 18, 2024}

\def\eqref#1{(\ref{#1})}

\newcommand{\Z}{{\Bbb Z}}
\newcommand{\C}{{\Bbb C}}

\newcommand{\R}{{\Bbb R}}

\def\1{\sqrt{-1}\:}

\makeatletter
\def\x@arrow{\DOTSB\Relbar}
\def\xlongequalsignfill@{\arrowfill@\x@arrow\Relbar\x@arrow}
\newcommand{\xlongequal}[2][]{%
        \ext@arrow 0099\xlongequalsignfill@{#1}{#2}}
\def\xlongrightarrowfill@{\arrowfill@\relbar\relbar\longrightarrow}
\newcommand{\xlongrightarrow}[2][]{%
        \ext@arrow 0099\xlongrightarrowfill@{#1}{#2}}
\makeatother

\renewcommand{\bar}{\overline}
\renewcommand{\phi}{\varphi}
\renewcommand{\epsilon}{\varepsilon}
\renewcommand{\geq}{\geqslant}
\renewcommand{\leq}{\leqslant}

\newcommand{\Id}{\operatorname{Id}}

\newcommand{\Kah}{\operatorname{Kah}}
\newcommand{\Pos}{\operatorname{Pos}}
\newcommand{\Aut}{\operatorname{Aut}}

\newcounter{Mycounter}[section]
\newcounter{lemma}[section]
\renewcommand{\thelemma}{{Lemma \thesection.\arabic{lemma}}}
\newcommand{\lemma}{%
    \setcounter{lemma}{\value{Mycounter}}
    \refstepcounter{lemma}
    \stepcounter{Mycounter}
    {\noindent \bf \thelemma:\ }}

\newcounter{claim}[section]
\renewcommand{\theclaim}{{Claim \thesection.\arabic{claim}}}
\newcommand{\claim}{%
    \setcounter{claim}{\value{Mycounter}}
    \refstepcounter{claim}
    \stepcounter{Mycounter}
    {\noindent \bf \theclaim:\ }}

\newcounter{sublemma}[section]
\newcounter{corollary}[section]
\renewcommand{\thecorollary}{{Corollary \thesection.\arabic{corollary}}}
\newcommand{\corollary}{%
    \setcounter{corollary}{\value{Mycounter}}
    \refstepcounter{corollary}
    \stepcounter{Mycounter}
    {\noindent \bf \thecorollary:\ }}

\newcounter{theorem}[section]
\renewcommand{\thetheorem}{{Theorem \thesection.\arabic{theorem}}}
\newcommand{\theorem}{%
    \setcounter{theorem}{\value{Mycounter}}
    \refstepcounter{theorem}
    \stepcounter{Mycounter}
    {\noindent \bf \thetheorem:\ }}

\newcounter{conjecture}[section]
\renewcommand{\theconjecture}{{Conjecture \thesection.\arabic{conjecture}}}
\newcommand{\conjecture}{%
    \setcounter{conjecture}{\value{Mycounter}}
    \refstepcounter{conjecture}
    \stepcounter{Mycounter}
    {\noindent \bf \theconjecture:\ }}

\newcounter{proposition}[section]
\renewcommand{\theproposition}
      {{Proposition \thesection.\arabic{proposition}}}
\newcommand{\proposition}{%
    \setcounter{proposition}{\value{Mycounter}}
    \refstepcounter{proposition}
    \stepcounter{Mycounter}
    {\noindent \bf \theproposition:\ }}

\newcounter{definition}[section]
\renewcommand{\thedefinition}{{Definition~\thesection.\arabic{definition}}}
\newcommand{\definition}{%
    \setcounter{definition}{\value{Mycounter}}
    \refstepcounter{definition}
    \stepcounter{Mycounter}
    {\noindent \bf \thedefinition:\ }}

\newcounter{example}[section]
\newcounter{remark}[section]
\renewcommand{\theremark}{{Remark \thesection.\arabic{remark}}}
\newcommand{\remark}{%
    \setcounter{remark}{\value{Mycounter}}
    \refstepcounter{remark}
    \stepcounter{Mycounter}
    {\noindent \bf \theremark:\ }}

\newcounter{problem}[section]
\newcounter{question}[section]
\newcommand{\proof}{{\bf Proof:\ }}

\makeatletter

\@addtoreset{equation}{section} \@addtoreset{footnote}{section} \makeatother

\def\blacksquare{\hbox{\vrule width 5pt height 5pt depth 0pt}}
\def\endproof{\blacksquare}

\begin{document}
\begin{center}
{\LARGE\bf
Special Hermitian structures on  suspensions\\[4mm]
}

Anna Fino, Gueo Grantcharov, Misha Verbitsky\footnote{Anna Fino
is partially supported by    Project PRIN 2022 \lq \lq Real and complex manifolds: Geometry and Holomorphic Dynamics”,  by  GNSAGA (Indam) and by a grant from the Simons Foundation (\#944448), Gueo
Grantcharov is partially supported by a grant from the Simons
Foundation (\#853269). Misha Verbitsky is partially supported by
CNPq - Process 310952/2021-2 and FAPERJ E-26/202.912/2018.}

\end{center}


{\small
\begin{minipage}[t]{0.85\linewidth}
{\bf Abstract}.   Motivated by  the construction based on
topological suspension of a family of compact
non-K\"ahler complex manifolds  with trivial canonical
bundle given  by L. Qin and  B. Wang in \cite{QW},
we   study  toric  suspensions  of balanced manifolds by
holomorphic automorphisms. In particular, we show that
toric  suspensions of Calabi-Yau manifolds
are balanced. We also prove that suspensions
associated with hyperbolic automorphisms
of hyperk\"ahler manifolds do not admit
 any pluriclosed, astheno-K\"ahler or
$p$-pluriclosed Hermitian metric.
Moreover, we consider natural extensions  for  hypercomplex
manifolds,  providing  some explicit examples of compact
holomorphic symplectic and hypercomplex non-K\"ahler
manifolds. We also show that a modified suspension
construction provides examples with pluriclosed metrics.
\end{minipage}
}


\tableofcontents


\section{Introduction}


Finding examples of special non-K\"ahler metrics  on
compact complex manifolds has become a question of
increasing interest in recent years. It is partly due to
developments in Physics related to Hull-Strominger system
(\cite{St}, \cite{Hu}) and generalized geometry
(\cite{Gua1}, \cite{Gua2}, \cite{Hi}). In \cite{QW} a new
class of examples of non-K\"ahler manifolds with trivial
canonical bundle and nice topological properties have been
introduced.  It is based on the topological suspension
construction.

Given in general  a  smooth  manifold $M$ and a
diffeomorphism $f$ of $M$, the mapping torus  (or
suspension) of $f$ is defined to be the quotient $M_f$ of
the product $M  \times \R$  by the $\Z$-action defined by
$$
(p, t) \rightarrow (f^{-n} (p), t + n).
$$
As a consequence $dt$ defines a nonsingular closed  $1$-form on  $M_f$  tangent to the fibration
$$
M_f \longrightarrow S^1 = \R/\Z.
$$
Moreover, the vector field $\frac{\partial }{\partial  t}$
on $M \times \R$  defines a vector field   on $M_f$, the
suspension of
the diffeomorphism $f$. There is a natural correspondence
between the orbits of  $f$  and the trajectories of  the
vector field.
Mapping tori have been used  in \cite{Li} to construct
examples of co-symplectic and co-K\"ahler manifolds.

The  suspension construction can be extended to complex
manifolds in the following way.
Given a  complex manifold $M$,  a set of   commuting
holomorphic automorphisms $f_j$, $j = 1, \ldots, 2k$,  of
$M$ and a lattice $\Lambda \subset \C^k$  of rank 2k,
generated by $\xi_1, \ldots , \xi_{2k}$, one can define an
action of $\Z^{2k} = \langle \xi_1, \ldots, \xi_{2k}
\rangle$  on $M  \times \C^k$ via
$\varphi_j(m, z) = (f_j(m), z + \xi_j )$. The quotient of $M \times  \C^k$ by the action  of $\Z^{2k}$ is called the toric suspension of $(M,f_1, \ldots, f_{2k}).$
In particular, if $f$ is an automorphism of a
complex manifold $M$ and $T^2 = \C/\Z^2$ an elliptic
curve, one  can construct   the complex suspension of $f$
as  the toric suspension  $S(f)$  of $M$ associated with
the pair $(f, {\text {Id}}_M).$ In the present paper  we  study the metric properties of the
constructed manifolds, like the existence of balanced
metrics, that is, Hermitian metrics    with co-closed
fundamental form.    We
also extend the construction to produce  hypercomplex
manifolds with special  metric properties.

Using a different construction related to automorphisms
of 3-dimensional Sasakian manifolds, we construct
suspensions admitting pluriclosed metrics, that is, the  Hermitian
metrics with $\partial \overline \partial$-closed
fundamental forms.

In Sections \ref{_hyperkahler_Section_} and \ref{section3} we present the necessary information
on  hyperk\"ahler manifolds and  toric suspension
construction.   In  Section \ref{Sect4}  we prove that the
complex toric  suspension of a balanced   manifold  $M$ by
two commuting holomorphic diffeomorphisms preserving a
volume form is balanced.  As a  corollary  we   show that
if  $M$ is a Calabi-Yau manifold and
$f$  is an automorphism  of $M$  preserving the
holomorphic volume form,  then the complex suspension
$S(f)$ has trivial canonical bundle and admits a balanced
metric.

In Section \ref {_k-plurikahler_} we show that the
balanced manifolds constructed using  any
hyperbolic automorphism of hyperk\"ahler
manifolds do not admit any $p$-pluriclosed and locally
conformally K\"ahler (LCK) metric.  In Section \ref{Sect5}  we
recover the construction in \cite{QW} as toric suspension
of a Kummer surface and we generalize it to suspension of
the  Hilbert scheme of points on Kummer surfaces.
In Section  \ref{_holo_symplectic-Section_} we discuss the
natural extensions of toric suspensions on hypercomplex
manifolds and their metric structures. As an application
we  construct  explicit examples of compact holomorphic
symplectic and hypercomplex non-K\"ahler manifolds. The
examples are in fact pseudo-hyperk\"ahler and admit
quaternionic balanced metric, but no  hyperk\"ahler  with torsion (HKT) metric.

Finally in   Section \ref{Sect7} we show how using
automorphisms of  Sasakian and K\"ahler manifolds  it is
also possible to construct suspensions admitting
pluriclosed metrics recovering a recent example
constructed in \cite{_Fino:Paradiso},  as a  compact
$3$-step solvmanifold.


\section{Hyperk\"ahler manifolds and their automorphisms}
\label{_hyperkahler_Section_}


Here we introduce the necessary background materials on hyperk\"ahler
geometry. We follow \cite{_AV:Aut_}, \cite{_Beauville_},
\cite{_Besse:Einst_Manifo_}, \cite{_Cantat:K3_},
\cite{_Kapovich:Kleinian_}.

\subsection{Hyperk\"ahler manifolds and the BBF form}

\definition
A {\bf hyperk\"ahler manifold}
is a compact, K\"ahler, holomorphically symplectic manifold.

\hfill

\definition
A hyperk\"ahler manifold $M$ is called to be
{\bf of maximal holonomy} (also: simple, or IHS)
if $\pi_1(M)=0$ and  $H^{2,0}(M)=\C$.

\hfill

\theorem (Bogomolov's decomposition \cite{_Bogomolov_})\\
 Any
hyperk\"ahler manifold admits a finite covering
which is a product of a torus and several
 hyperk\"ahler manifolds of maximal holonomy.

\hfill

\remark
From now on  all hyperk\"ahler manifolds
are assumed to be of maximal holonomy.

\hfill

\theorem (Fujiki \cite{_Fujiki_})\\
Let  $M$ be a hyperk\"ahler manifold, $\eta\in H^2(M, \Z)$, and  $n = \frac{\dim_{\C} M}{2}.$  Then $\int_M \eta^{2n}=cq(\eta,\eta)^n$,
where $q$ is a primitive integer non-degenerate  quadratic form  on $H^2(M,\Z)$,
and $c>0$  is a rational number depending only on $M$.

\hfill

\definition
This  primitive integral quadratic form $q$ on  $H^2 (M, \Z)$  is called
{\bf Bogomolov-Beauville-Fujiki form}, or {\bf BBF form}.
It is defined
by the Fujiki's relation uniquely, up to a sign. The sign is determined
from the following formula (Bogomolov, Beauville, see \cite{_Beauville_})
\begin{align*}  \lambda q(\eta,\eta) &=
   \int_M \eta\wedge\eta  \wedge \Omega^{n-1}
   \wedge \bar \Omega^{n-1} \\
 &-\frac {n-1}{n}\left(\int_M \eta \wedge \Omega^{n-1}\wedge \bar
   \Omega^{n}\right) \left(\int_M \eta \wedge \Omega^{n}\wedge \bar \Omega^{n-1}\right),
\end{align*}
where $\Omega$ is the holomorphic symplectic form on $M$ and
$\lambda>0$.

\hfill

\remark
The BBF form $q$ has signature $(3,b_2-3)$ when extended on $H^2 (M,  \R)$.
It is negative definite on primitive forms, and positive
definite on $\langle \Omega, \bar \Omega, \omega\rangle$,
 where $\omega$ is a K\"ahler form. On (1,1)-forms $\eta$ it
can be written as $q(\eta, \eta) =  c \int_M
\eta\wedge\eta  \wedge \Omega^{n-1} \wedge \bar \Omega^{n-1}$,  where $c$ is a constant.

\subsection{Classification of automorphisms of hyperk\"ahler manifolds}

\remark The indefinite orthogonal group $O(m,n),$ $m, n>0$,  is the Lie group of all linear transformations of an $l$-dimensional real vector space that leave invariant a nondegenerate, symmetric bilinear form  $q$ of signature $(m, n)$, where $l = m + n$.   $O(m, n)$ has 4 connected components.
We denote the connected component of 1 by $SO^+(m,n)$.
We call a vector $v$ {\bf positive} if its square is positive, i.e. if $q (v, v) >0$.

\hfill

\definition
Let $V$ be a real vector space of dimension $n + 1$   with a  quadratic form $q$ of signature
$(1,n)$, $\Pos(V)=\{x\in V\ \ |\ \ q(x,x)>0\}$
its {\bf  positive cone},  and ${\Bbb P}^+ V$  be the projectivization
of $\Pos(V)$.
Denote by $g$ any $SO(V)$-invariant Riemannian structure on
 ${\Bbb P}^+ V$ (it is easy to see that $g$ is
unique up to a constant multiplier). Then $({\Bbb P}^+ V, g)$ is called
{\bf hyperbolic space}, and $SO^+(V)$
{\bf the group of oriented  hyperbolic isometries}.

\hfill

\theorem \label{th2.10}
Let $n>0$, and $\alpha \in SO^+(1,n)$ be  an isometry
acting on $V$. Then one and only one of the following  three cases occurs
\begin{description}
\item[(i)] $\alpha$ has an eigenvector $x$ with $q(x)>0$
($\alpha$ is {\bf ``elliptic isometry''});
\item[(ii)]
 $\alpha$ has  two  eigenvectors $x$ and $y$ such that  $q(x,x)= q(y,y) =0$
and  real eigenvalues $\lambda_x$ and $\lambda_y= \lambda_x^{-1}$  satisfying $|\lambda_x|>1$ and all other eigenvalues have absolute value one
($\alpha$ is {\bf ``hyperbolic isometry''}, or {\bf loxodromic isometry});
\item[(iii)] $\alpha$ has a unique (up to a  constant) eigenvector $x$ with $q(x,x)=0$ with eigenvalue 1, and
no fixed points on  ${\Bbb P}^+ V$
($\alpha$ is {\bf ``parabolic isometry''}).
\end{description}

For a proof see for instance \cite{_Kapovich:Kleinian_} or \cite[Chapter 5]{Gallier:diffgeom_and_Lie}.

\hfill

\remark
 All eigenvalues
of elliptic and parabolic
isometries have absolute value 1. Hyperbolic and elliptic isometries
are semisimple (that is, diagonalizable over $\C$).

\hfill

\definition
Notice that any complex automorphism of a hyperk\"ahler manifold
acts by isometry on the space $H^{1,1}(M,\R)$ with the BBF metric
which has signature $(1, b_2-3)$.
A complex automorphism $f$ of a hyperk\"ahler manifold $M$
is called {\bf elliptic (parabolic, hyperbolic)}
if the  induced action  $f^*$ of $f$ is elliptic (parabolic, hyperbolic) on $H^{1,1}(M,\R)$.

\hfill

Further on we shall need the following lemma.

\hfill

\lemma\label{_nega_defi_on_inva_Lemma_}
Let $M$ be a hyperk\"ahler manifold,
$f:\; M \longrightarrow M$ a hyperbolic automorphism,
and $\eta\in H^{1,1}(M, \R)$ a non-zero $f^*$-invariant class.
Then $q(\eta, \eta) < 0$.

\hfill

{\bf Proof:}  Let $v_+, v_-$ be eigenvectors of $f^*$
with the real eigenvalues $\lambda >1$ and $\lambda^{-1}$.
Then any invariant vector of $f^*$ belongs to
$\langle v_+, v_-\rangle^{\perp}$. However, $q$
is negative definite on the space spanned by the other eigenvectors because
signature of $q$ on $H^{1,1}(M,\R)$ is $(1, b_2-3)$.
\endproof

\newpage


\section{Toric suspensions} \label{section3}


\subsection{Toric suspensions: definition and basic properties}

\definition
Let $M$ be a complex manifold, and $f_1, ..., f_{2k}\in \Aut(M)$
a set of  commuting  holomorphic automorphisms of $M$.
Let $\Lambda\subset \C^k$ be a lattice of rank $2k$,
generated by $\xi_1, ..., \xi_{2k}$.
Define an action of $\Z^{2k}= \langle \phi_1, ..., \phi_{2k} \rangle$
on $M\times \C^k$ via $\phi_j(m, z)= (f_j(m), z + \xi_j )$.
In other words, $\Z^{2k}$ acts on $\C^k$ as a shift by
the corresponding element of $\Lambda$ and on
$M$ as an automorphism obtained as an appropriate
product of $f_i$. The  quotient $(M\times \C^k)/\Z^{2k}$
is called {\bf the toric suspension} of $(M, f_1, ..., f_{2k})$.

\hfill

\remark
The toric suspension is clearly complex analytic, holomorphically
fibered over the torus $\C^{k}/\Lambda$, but not necessarily K\"ahler.

\hfill

\theorem\label{_suspension_non-Kahler_Theorem_}
Let $S(M, f_1, \ldots, f_{2k})$ be
a toric suspension, with $M$ a compact K\"ahler
manifold. Then $S(M, f_1, \dots, f_{2k})$ is K\"ahler
if and only if there is a K\"ahler class $[\omega]\in H^{1,1}(M)$
such that $f_i^*([\omega])=[\omega]$.

\hfill

{\bf Proof:}  See  the proof of Theorem 3.4.1 in the paper \cite{_Manjarin_}.
\endproof


\subsection{Hyperbolic suspensions}

The following definition is motivated by
the classification of the automorphism groups of hyperbolic
manifolds, such as a K3 surface.

\hfill

\definition \label{def3.4}
Let $f:\; M \longrightarrow M$ be an automorphism of a compact
complex manifold of K\"ahler type (i.e. admitting a K\"ahler metric).
We say that $f$ is {\bf  a hyperbolic
automorphism} if the  induced action  of $f$ on $H^{1,1}(M, \R)$
has a unique (up to a constant)  eigenvector $\eta$ with  eigenvalue 
 $\lambda>1$.

\hfill

We list some immediate properties of hyperbolic automorphisms.

\hfill

\proposition\label{_hype_basic_Proposition_}
Let $f:\; M \longrightarrow M$ be a hyperbolic automorphism of a compact
complex manifold of K\"ahler type, and $\eta\in H^{1,1}(M)$
an eigenvector with an eigenvalue $f^*\eta= \lambda \eta$
such that $\lambda>1$. Denote by $\Kah(M)\subset H^{1,1}(M, \R)$
the K\"ahler cone of $M$. Then
\begin{description}
\item[(i)] $\eta$ belongs to the closure of the K\"ahler cone.
\item[(ii)] $\int_M \eta^{n}=0$, where $n=\dim_\C M$.
In particular, $\eta\notin\Kah(M)$.
\item[(iii)] The action of $f$ on $\Kah(M)$ has no fixed points.
\end{description}
{\bf Proof:}
Let $S\subset H^{1,1}(M,\R)$
be the sum of all eigenspaces of $f$ on $H^{1,1}(M, \R)$
with eigenvalues  different from   $\lambda$. Since $\lambda$ is the biggest eigenvalue,
for any $v\in H^{1,1}(M,\R)\backslash S$, one has
$\lim_i \frac{(f^*)^i(v)}{\lambda^i}=c \eta$  where $c$ is  a nonzero constant. Since $\Kah(M)$
is open, this is also true for general K\"ahler class $\omega$.
We obtained $\eta$ as a limit of K\"ahler forms. This proves (i).

To see (ii), we notice that $\int_M \eta^{n}= \int_M
f^*(\eta)^{n}= \lambda^{n} \int_M \eta^{n}$.

To obtain (iii), assume that $f$ fixes a K\"ahler class
$\omega$ on $M$.
Then $f$   is an elliptic isometry on $H^{1,1}(M,\R)$, but  by  \ref{th2.10}
$f$ can not be hyperbolic, giving a contradiction.
\endproof

\hfill

\remark
Since a  hyperbolic automorphism of a hyperk\"ahler manifold preserves
its K\"ahler cone, and the eigenvector $x$ with $|\lambda_x|>1$
sits on the boundary of the K\"ahler cone (\ref{_hype_basic_Proposition_}),
the number $\lambda_x$ is positive.

\hfill

\definition
Let $f:\; M \longrightarrow M$ be an automorphism
of a  compact complex manifold of K\"ahler type,  and $T^2= \C /\Z^2$
an elliptic curve. Consider a toric suspension $S(f)$ of $M$
associated with the pair $(f, \Id_M)$. This manifold is called
{\bf a complex suspension of $f$}. We call $S(f)$
{\bf a hyperbolic suspension} if $f$ is hyperbolic.

\hfill

\remark\label{_suspension_explicit_Remark_}
  The toric suspension $S(f)$ of $M$
associated with the pair $(f, \Id_M)$ can be viewed as the product  manifold  $M_f  \times S^1$, where $M_f$ is the mapping torus of $M$ by $f$ obtained as the quotient of $M \times \R$ by the $\Z$-action
 $$
(p, t)   \rightarrow  (f^{-n}(p), t + n).
$$
If $(t,s)$ are local coordinates on $\R \times S^1$, then
$\frac{\partial}{\partial t}$ on $M \times \R$ defines a
vector field $X_f$ on $S(f)$ called the {\bf suspension vector field} of
$f$ (see \cite{_Goldman_}). Note that the vector field
$X_f - i \frac{ \partial}{\partial s}$ on $S(f)$ is
holomorphic.  Moreover the vector fields $X_f,
\frac{\partial}{\partial s}$ provide a natural splitting
$TS(f) = T_{vert}S(f) \oplus \pi^*TE$, which defines a flat
Ehresmann connection on $S(f)$, which we call the {\bf
  standard connection}. We will denote by $\theta$ the
associated connection $1$-form such $\theta  + \1 ds$ is a
$(1,0)$-form with respect to the complex structure on
$S(f)$.

\hfill

\remark \label{rem3.9}
By \ref{_hype_basic_Proposition_} (iii)
and \ref{_suspension_non-Kahler_Theorem_},
a hyperbolic suspension is never K\"ahler.

\section{Balanced metrics on Calabi-Yau
suspensions}\label{Sect4}

Balanced  metrics were introduced  in \cite{_Michelson_}. For  further   properties and examples see e.g.   \cite{AB}, \cite{AB2} and \cite{FLY12}.

\hfill

\definition
Let $(M,I, h)$ be a complex Hermitian manifold, $\dim_\C M=n$, and
$\omega$ the fundamental $(1,1)$-form associated to $h$.
We say that $h$ is {\bf balanced} if $\omega^{n-1}$ is closed.

\hfill

The main result of the present Section is the following theorem.

\hfill

\theorem\label{_CY_susp_balanced_main_Theorem_}
Let $M$ be a balanced compact manifold  of complex dimension $n$ and $f_1, f_{2} \in {\mbox{Aut}}(M)$  two commuting   holomorphic automorphisms preserving a volume form $V$.
Denote by $\pi:\; S\to E$ the corresponding suspension over
an elliptic curve $E$. Assume that $M$ is
balanced. Then $S$ is also balanced.

\hfill

{\bf Proof.}
Let $\omega_E$ be a K\"ahler form on $E$.
Recall that a smooth fibration $\pi:\; S\to E$ over an elliptic curve
is called {\bf essential} (\cite{_Michelson_})
if $\pi^*(\omega_E)$ is not Aeppli exact, i.e. $\pi^*(\omega_E)$ cannot be  equal to $ \overline \partial \alpha + \partial \overline \alpha$, for any $(1,0)$-form $\alpha$.
Michelsohn (\cite {_Michelson_}) proves that the total space $S$
of an essential fibration with balanced fibers over a complex curve  is balanced.
 To prove   \ref{_CY_susp_balanced_main_Theorem_}
it remains only to show that $\pi^*(\omega_E)$ is not Aeppli
exact. 

Since $V$ is $f_j$-invariant, $ j =1,2$, we  may extend $V$ to a form $V_h$ on $S$ vanishing on horizontal
vector fields of this Ehresmann connection.
Then the  form  $V_h$  is  of type $(n, n)$, positive and closed.
 Since $V_h$ vanishes on any horizontal vector, the form $\pi^*(\omega_E)\wedge V_h$ is of maximal degree and positive, so $\int_S \pi^* \omega_E \wedge V_h>0$. To prove  the theorem by contradiction assume that $\pi^*(\omega_E)$  is Aeppli exact. However by  Stokes Theorem we would  have $\int_S \pi^* \omega_E \wedge V_h=0$, which is impossible.
\endproof


\section{Hyperbolic holomorphically symplectic suspensions} \label{_k-plurikahler_}


\subsection{Hyperbolic holomorphically symplectic suspensions}



\definition
Let $M$ be a hyperk\"ahler manifold and
$f:\; M \longrightarrow M$ a hyperbolic automorphism (as in  \ref{def3.4})
preserving the holomorphic symplectic form.
Denote by $S$ the corresponding hyperbolic suspension,
fibered over $T^2$ with fiber $M$.
Then $S$ is called {\bf a  hyperbolic holomorphically symplectic
suspension}.

\hfill

Similarly, if $M$ is a Calabi-Yau manifold and  $f$ is a hyperbolic  automorphism of  $M$  preserving the complex  holomorphic volume form, we will call $S(f)$ a {\bf Calabi-Yau hyperbolic suspension}.

\hfill

\proposition  \label{_HS_susp_balanced_main_Theorem_} Let S be a hyperbolic holomorphically symplectic suspension or a Calabi-Yau hyperbolic suspension. Then $S$ is balanced and non-K\"ahler Calabi-Yau.

\smallskip

\proof  In both cases there exists an invariant non-vanishing  holomorphic section   $\Theta$  of the
canonical bundle of $M$.
Therefore, $V:=\Theta \wedge \bar \Theta$ is a
$f$-invariant volume on $M$.  By
\ref{_CY_susp_balanced_main_Theorem_} $S$ is balanced.  By
\ref{rem3.9} $S$ is non-K\"ahler. Moreover,  the form
$(\theta  + \1 ds) \wedge \Theta$
(\ref{_suspension_explicit_Remark_})
is a non-vanishing holomorphic section of the  canonical bundle of $S$.
\endproof

\hfill


\subsection{Balanced, pluriclosed and LCK Hermitian metrics}

The study of special Hermitian metrics   posed also the question of compatibility between different structures of non-K\"ahler type.  We recall the  conjecture in \cite{FV16} according to which a compact complex manifold admitting both a pluriclosed, i.e.  whose  Hermitian form $\omega$  satisfies $dd^c \omega=0$,  and a balanced metric is K\"ahler. This has been already proven for  specific cases in the papers \cite{Ve14, C14, FV16, FLY12,Otiman, FP, GiPo,_Freibert-Swann_}.
A similar question  was posed in   \cite{STW}  (see also \cite{Fe}) for  a compact complex manifold  of complex dimension $n$  admitting a balanced metric and an
astheno-K\"ahler metric, i.e.   whose Hermitian  form satisfies $d d^c \omega^{n -2}=0$. A negative answer to this question  was given  in \cite{FGV, LU}. For conjectures related to  the existence of locally conformally K\"ahler metrics - the ones that satisfies $d\omega=\theta \wedge \omega$, see the book  \cite{OV}.

\hfill

Based on the previous discussion one can formulate the following general conjecture:

\hfill

\conjecture
Let $X$ be a compact complex manifold, $n: = \dim_\C X > 2$.
Assume that two of the following assumptions occur.
\begin{description}
\item[(i)] $X$ admits a Hermitian form $\omega$ which
is locally conformally K\"ahler, that is, satisfies $d\omega=\theta \wedge \omega$.
\item[(ii)] $X$ admits a Hermitian form $\omega$ which
is balanced.
\item[(iii)]  $X$ admits a Hermitian form $\omega$ which
is $p$-pluriclosed, that is, satisfies $dd^c(\omega^p)=0$,  for  $p= 1, 2, \ldots, n-3$ if $n >3$ or  for $p= 1$ if $n =3$.
\end{description}
Then $X$ admits a K\"ahler structure.

\hfill

In this section, we prove this conjecture
when $X$ is a suspension of a hyperk\"ahler manifold $M$
associated with a hyperbolic automorphism of $M$. The non-existence of locally conformally  K\"ahler metrics locally but not globally conformal on these examples follows  from Proposition 37.8 in \cite{OV}.

\subsection{Strongly positive and weakly positive $(p,p)$-currents}

Here we recall that a $(p,p)$-{\em current} on a complex manifold $X$ is an element of the Frechet space dual to the space of $(n-p,n-p)$ complex forms $\Lambda^{n-p,n-p}(X)$. In the compact case, the space of $(p,p)$-currents can be identified with the space of $(p,p)$-forms with distribution coefficients and the duality is given by integration. So for any $(p,p)$-current $T$ and a form $\alpha$
of type $(n-p,n-p)$ we have
$$
\langle T,\alpha\rangle = \int_X T\wedge \alpha.
$$
The operators $d$ and $d^c$ can be extended to $(p,p)$-currents by using the duality induced by the
integration, i.e.,   $dT$ and $d^cT$ are respectively defined via the relations
$$
\langle dT,\beta\rangle = - \int_X T\wedge d\beta, \quad \langle d^cT,\beta\rangle = -\int_X T\wedge d^c\beta.
$$

\hfill

We recall now the definition of a positive $(p, p)$-form (see e.g. \cite[Chapter 3]{_Demailly:agbook_}).

\hfill

\definition A {\bf weakly positive}   ({\bf  strictly weakly positive})  $(p, p)$-form on a complex manifold  $X$  is a real $(p, p)$-form $\eta$  such that
 for any complex subspace  $V  \subset TM,$ $\dim_C V = p,$ the restriction  $\eta  \mid_V$  is a non-negative volume form (positive volume form).  Weakly positive condition is equivalent to
 $$
 i^p \eta (v_1, \overline v_1, v_2, \overline v_2, \ldots v_p, \overline v_p) \geq 0,
 $$
 for every tangent vectors $v_1, \ldots, v_p \in T^{1,0}_x X$. A real $(p,p)$-form $\eta$ is called {\bf strongly positive}  ({\bf strictly strongly positive}) if it can be locally expressed as a sum
$$
\eta = i^p  \sum_{j_1, \ldots, j_p} \eta_{j_1 \ldots j_p} \xi_{j_1} \wedge \overline \xi_{j_1} \wedge \ldots \wedge \xi_{j_p } \wedge \xi_{ \overline j_p},
$$
 running over the  set of $p$-tuples $\xi_{j_1}, \xi_{j_2}, \ldots, \xi_{j_p}$ of $(1,0)$-forms, with $\eta_{j_1 \ldots j_p} \geq 0$ ($\eta_{j_1 \ldots j_p} > 0$).
\hfill

All strongly positive forms are also weakly positive. The strongly positive and the weakly positive forms form closed, convex cones in the space of real $(p, p)$-forms,  see for instance \cite[Chapter 3]{_Demailly:agbook_}. These two cones are dual with respect to the  natural  pairing
$$
\Lambda_x^{p,p} (X, \R) \times \Lambda_x^{n-p,n-p} (X, \R)  \rightarrow  \R.
$$ For $(1,1)$-forms and $(n-1, n-1)$-forms, the strong positivity is equivalent to weak positivity. Finally, a  product of a weakly positive form and a strongly positive one is always weakly positive (however, a  product of two weakly positive forms may be not weakly positive).  A product of strongly positive forms  is still strongly positive.

\hfill

A strongly/weakly positive $(p,p)$-current is a current taking non-negative values on weakly/strongly positive compactly supported $(n - p, n - p)$-forms.

\hfill

\definition A $(p,p)$-current $T$ is called {\bf weakly positive}
if
$$
i^{n-p}\int_X T\wedge\alpha_1\wedge\overline{\alpha_1}\wedge...\alpha_{n-p}\wedge\overline{\alpha}_{n-p} \geq 0,
$$
 for every (1,0)-forms $\alpha_1,...\alpha_{n-p}$ with inequality being strict for at least one choice of $\alpha_i$'s. The current $T$ is called {\bf strongly  positive} if the inequality is strict for every non-zero $\alpha_1\wedge\overline{\alpha_1}\wedge...\alpha_{n-p}\wedge\overline{\alpha}_{n-p}$.

\hfill

\definition
 A $(p,p)$-current $T$ is called {\bf  strictly strongly  positive}   (resp. {\bf  strictly weakly   positive}) if $T > \epsilon \omega$ for a strictly  strongly positive  (resp. strictly weakly   positive)  $(p,p)$-form  $\omega$ and a positive number $\epsilon$.

\hfill

\claim
The cone of strongly positive $(p,p)$-currents is
dual to the cone of strictly weakly positive $(p,p)$-forms,
the cone of weakly positive $(p,p)$-currents is
dual to the cone of strictly strongly positive $(p,p)$-forms.

%
\hfill
%
\endproof

The main result of this section is the following

\theorem\label{_dd^c_exact_suspension_Theorem_}
Let $f\in \Aut(M)$ be a hyperbolic automorphism of a
hyperk\"ahler manifold, and denote by
$\pi:  S \rightarrow E$ the suspension $S(f)$ of $(M, f)$.
Then $S$ admits a $dd^c$-exact, strongly positive
$(p,p)$-current $\beta$ for any $p= 2, 3, \ldots, n-1$, where  $n:=\dim_\C M$.

\hfill

We prove this theorem in Subsection \ref{_CDS_currents_Subsection_}.
\ref{_dd^c_exact_suspension_Theorem_} immediately implies the following.

\hfill

\corollary \label{_Corollary5.7_}
Let $S$ be a hyperbolic suspension over a hyperk\"ahler
manifold $M$, as in \ref{_dd^c_exact_suspension_Theorem_}.
Then $S$ does not admit a $dd^c$-closed strictly weakly positive
$(n-p+1, n-p+1)$-form $U$ for $p= 2,3,\ldots,n-1$. In particular, $S$ is not $k$-pluriclosed
for any $k=1, 2, ..., n-1$.

\proof Let $\beta= dd^c \alpha$ be a current introduced in
\ref{_dd^c_exact_suspension_Theorem_}.
If $U$ is $dd^c$-closed strictly weakly positive $(n-p+1, n-p+1)$-form $U$,
we have $0< \int_M U\wedge \beta =\int_M dd^c U\wedge \alpha =0$,
which is impossible.
\endproof

\subsection{Hyperbolic automorphisms and Cantat-Dingh-Sibony currents}
\label{_CDS_currents_Subsection_}

Let $f$ be a hyperbolic automorphism of a hyperk\"ahler manifold $M$,
$\dim_\C M=n$, and $p= 1, 2, ..., n-1$, and denote by $\lambda$
its   unique eigenvalue which satisfies $|\lambda| >1$. 
 
 Recall that
{\bf the mass} of a positive $(p,p)$-current $v$ 
on a K\"ahler manifold $M$ is $\int_M v \wedge \omega^{n-p}$.
Since $f$ preserves the K\"ahler cone, it preserves the positive
cone of $M$, hence $\lambda >1$.
The action of $f^*$ on
$H^{2p}(M)$ has $\lambda^p$ as the maximal eigenvalue
(\cite{_BKLV_}), hence the mass of $\frac 1 {\lambda^{pk}}
(f^*)^k\omega^p$ is bounded.
Moreover,  the set of positive $(p,p)$-currents of bounded mass is compact
(\cite[Chapter 3]{_Demailly:agbook_}).

 Therefore
the sequence $\{ \frac 1 {\lambda^{pk}} (f^*)^k\omega^p \}_{k =1, \ldots, \infty}$
has a limit point.
The eigenspace corresponding to $\lambda^p$ in
$H^{p,p}(M)$ has multiplicity 1, as shown in \cite{_BKLV_}.
By \cite[Theorem 4.3.1]{DinhSibony10},
the limit  of a subsequence $\lim_k\frac 1 {\lambda^{pk}} (f^*)^k\omega^p$
is a unique positive $(p,p)$-current $\sigma$
which satisfies $f^*\sigma = \lambda \sigma$.
We call it {\bf the Cantat-Dingh-Sibony current}
(Cantat prove this result for  (1,1)-currents
on a K3 surface \cite{_Cantat:K3_}, and Dingh-Sibony for
all dimensions).

Using the decomposition
$TS = T_{{vert}} S \oplus \pi^* TE$ induced  by the
flat Ehresmann connection on $S$, we can
consider the bundle ${\Bbb D}:=D^{p,p}_\pi(S)$
of fiberwise currents as a local system
on $E$; the monodromy of this local system
is given by the map $v \mapsto f^* v$.
Identifying local systems and flat bundles,
we can consider ${\Bbb D}$ as a bundle
with flat connection.

Consider a real line bundle $L\subset {\Bbb D}$
spanned by the Cantat-Dingh-Sibony
current $\sigma$. This line bundle is preserved
by the natural flat connection on ${\Bbb D}$,
and its monodromy map is multiplication by $\lambda$.
Choose a trivialization of $L$ such that
the corresponding connection $1$-form $\theta$ satisfies $d \theta =0$ and $d (I \theta) =0$.
Using the decomposition $TS = T_{vert}  S \oplus \pi^* TE$,
we can embed the sections of ${\Bbb D}$ into
the space of $(p,p)$-currents on $S$. Let
$\alpha$ be the current on $S$ associated
with the section of $L\subset {\Bbb D}$ constructed above.
Since  $\sigma$ is a limit of closed currents,   $d\sigma=0$ and we have $d\alpha = \alpha \wedge\theta$, and
$\beta:=dd^c \alpha = \alpha \wedge\theta\wedge I \theta$. The current $\beta$
is $dd^c$-exact. Since $\beta$ is a limit
of the wedge power  of  strongly positive (1,1)-forms,
it is strongly positive. This proves
\ref{_dd^c_exact_suspension_Theorem_}.


\section{Examples of suspensions of hyperk\"ahler manifolds} \label{Sect5}

We briefly recall the examples of suspensions of Kummer K3 surfaces from  \cite{QW} first.

Take the complex  $2$-torus   ${\mathbb T}$ given by the  quotient of ${\mathbb C}^2$ by the standard  lattice generated by the unit vectors $(1,0),(i,0),(0,1),(0,i)$. Consider  the involution  of ${\mathbb C}^2$ given by multiplication by $-1$, i.e. $(z_1,z_2)\rightarrow(-z_1,-z_2)$. The involution descends to an involution $\sigma$ of the torus ${\mathbb T}$ with 16 fixed points $p_1, \ldots, p_{16}$.  The quotient  space ${\mathbb T}/<1, \sigma>$  has   16 double points. The singularities can be resolved by  blowing    the singularities  up, yieldings a smooth compact surface containing 16  mutually disjoint smooth rational curves $C_j$.   This is the  Kummer surface ${\rm {Km}}$ associated to $\mathbb T$. There is an alternative description of the Kummer surface. Let $X$ denote the surface obtaining by blowing up ${\mathbb T}$ at each of the points $p_1, \ldots, p_{16}$. Let $E_j \cong {\mathbb P}^1$ be the exceptional divisor over $p_j$. The involution $\sigma$ of ${\mathbb T}$ lifts to an involution  $\tau$ of $X$ with the fixed set $E = E_1 \cup \ldots \cup E_{16}$. The  eingevalues of the differential of $\tau$ at every points  of $E$ are $\pm 1$. So the quotient $X/<1, \tau>$ is smooth and contains 16 rational $(-2)$-curves $C_j \cong {\mathbb P}^1$, the images of the rational $(-1)$-curves $E_j$ in $X$. The quotient is a  Kummer surface ${\rm {Km}}$. Let $\hat{\mathbb C^2}$ be the surface obtained by blowing up ${\mathbb C}^2$ at every point of the discrete set $\pi^{-1} (\{ p_1, \ldots, p_{16} \})$, where $\pi:  {\mathbb C}^2 \rightarrow {\mathbb T}$ is the quotient map, we have the following diagram
\[
\begin{diagram}
\node{\hat{\mathbb C}^2} \arrow{e,t}{}  \arrow{s,l}{} \node{X}   \arrow{s,l}{}  \arrow{e,t}{}   \arrow{s,l}{}  \node{{\rm {Km}}}   \arrow{s,l}{\tilde \pi}\\
\node{{\mathbb C}^2}  \arrow{e,t}{\pi}   \node{\mathbb T} \arrow{e,t}{}   \node{{\mathbb T}/<1, \sigma>}
\end{diagram}
\]

By the Lefschetz Theorem on $(1,1)$-forms   we have  that the Picard group of ${\rm {Km}}$ is isomorphic to $H^2 ({\rm {Km}}, \Z) \cap H^{1,1}({\rm {Km}})$, so the rank of the Picard group of  ${\rm {Km}}$ is  $20$. Moreover, the Picard group of ${\rm {Km}}$  is generated by the 16 exceptional divisors $E_i$ and by the pull-back  by $\tilde \pi$  of divisors on ${\mathbb T}/<1, \sigma>$.

The canonical $(2,0)$-form $dz_1 \wedge dz_2$ on $\C^2$ induces a nowhere vaninshing $(2,0)$-form on ${\mathbb T}$. Therefore, the pullback of this form on $X$ induces a holomorphic $(2,0)$-form on the Kummer surface.

 Let $A\in SL(2, \Z+\sqrt{-1} \Z)$ be a matrix with $|tr(A)|>2$, so that it is diagonalizable with eigenvalues $\lambda,\lambda^{-1}$. Let   $dv_1, dv_2$  be respectively the associated  eigenvectors of the induced map on $H^1({\mathbb T},\C)\cong \Lambda^1(\R^4)$. Denote by $A$ also the induced map on $\Lambda^k(\R^4)$. Then $A$ preserves the holomorphic (2,0)-form $dv_1\wedge dv_2$ on $\mathbb{T}$ and the divisor $D=\sum_{i=1}^{16} E_i$. So it defines a holomorphic transformation $\phi_A$ on ${\rm {Km}}$ preserving the induced holomorphic $(2,0)$-form. In particular the $\mathbb{Z}$-action on $\mathbb{T}\times\mathbb{R}\times S^1$ generated by
 \begin{equation}\label{eq1}
 f: (p,x,y)\rightarrow (A(p),x+1,y)
 \end{equation}
 extends to an action on ${\rm {Km}}\times \mathbb{R}\times S^1$, generated by the hyperbolic automorphism $f$. The quotient is a  compact complex manifold $S({\rm {Km},f})$ with trivial canonical bundle and satisfying the hard Lefschetz property,  such that its real homotopy type is formal as shown in \cite{QW}.

 In a similar way we can construct hyperbolic automorphisms preserving the holomorphic symplectic form  on higher-dimensional hyperk\"ahler manifolds arising as Hilbert scheme of points on ${\rm {Km}}$. More precisely, $f: {\rm {Km}} \rightarrow {\rm {Km}}$ extends to $ f^{[n]} : {\rm {Km}}^{[n]}\rightarrow {\rm {Km}}^{[n]}$ on the Hilbert scheme of order $n$ of ${\rm {Km}}$ in a natural way: to a zero-dimensional subscheme $Z\subset {\rm {Km}}$ we assign $f(Z)$. According to Beauville (and PhD thesis by P. Beri   \cite{Beri}), $f^{[n]}$ preserves the holomorphic symplectic form if and only if $f$ does.
Now we can construct the suspension 
$S({\rm {Km}}^{[n]},  f^{[n]})$  using $f^{[n]}$ and obtain:

\hfill

\begin{theorem}
The space $S({\rm {Km}}^{[n]},  f^{[n]})$ for $n\geq 1$ is a non-K\"ahler compact complex manifold with trivial canonical bundle which admits a balanced metric and  it  is not $k$-pluriclosed for any $k = 1, 2, \ldots, 2n - 1$.\end{theorem}

\hfill

\proof The fact that it is balanced follows from \ref{_HS_susp_balanced_main_Theorem_} and \ref{_Corollary5.7_}.
\endproof

\hfill

The metric in the examples above is not explicit. But if we consider the suspension over the  real 4-torus $\mathbb{T}^4$ we can define such metric explicitly.  Denote by $v_1$ and $v_2$ the eigenforms of the map  $A$ on $H^1 (\mathbb{T}^4, \mathbb C)$ induced by the matrix $A$ as above and  by $x$ and $y$  respectively coordinates on $\R$ and $S^1$.  Then $A(dv_1\wedge \overline{dv_1})=|\lambda|^2dv_1\wedge \overline{dv_1}$    and $A(dv_2\wedge \overline{dv_2})=|\lambda |^{-2}dv_2\wedge \overline{dv_2}$. Consider the differential  forms on $\mathbb{T}^4 \times \R\times S^1$ given by $\alpha_1=|\lambda|^{-2x}dv_1\wedge d\overline v_1$ and $\alpha_2=|\lambda|^{2x} dv_2\wedge d\overline v_2$.  The forms $\alpha_1$ and $\alpha_2$ are  invariant under the action in (\ref{eq1}). Moreover  $\alpha_1+\alpha_2$ descends to a weakly  positive definite $(2,2)$-form on the suspension $S(\mathbb{T}^4, f)$ of the 4-torus defined by this action.  By the observation of Michelson \cite{_Michelson_} $S(\mathbb{T}^4, f)$ admits a balanced metric.  We can directly check that
$$ \omega =  |\lambda|^{-2x}dv_1\wedge d\overline{v_1} + |\lambda|^{2x}dv_2\wedge d\overline{v_2} + dx\wedge dy$$
is invariant and satisfies $d\omega^2 = 0$. Hence it defines a balanced metric.
\smallskip

\remark We restrict ourselves here to the more explicitly described examples, but many of the known compact hyperk\"ahler manifolds admit hyperbolic automorphisms. We expect that the topological properties of $S({\rm {Km},f})$ from \cite{QW} are also valid for $S({\rm {Km}}^{[n]},  f^{[n]})$. Note that   the manifold $A({\mathbb{T}}^4)$  can be also described as  the  almost abelian solvmanifold $M^6 (c)$ in \cite{AFLM} (see also Section 3 in \cite{FMS}). By Theorem 4.1  in \cite{FP} the  associated  almost abelian  Lie algebra, which is isomorphic to $\frak b_6$ in the notation of  \cite{AFLM},   admits  a balanced metric.  \hfill


\section{Holomorphically symplectic
and hypercomplex structures on toric suspensions with 4-dimensional base}
\label{_holo_symplectic-Section_}


In this section we show how the toric suspensions could be used to construct examples of compact holomorphic symplectic and hypercomplex non-K\"ahler manifolds. The examples are in fact pseudo-hyperk\"ahler. We also discuss their metric structure.

\subsection{General construction and example}

We consider a toric suspension of $M$, where $M$ is  a real  8-torus, over a  real 4-torus base.
More precisely we consider $S = S(T^8, f, id, id, id)$, where $f$ is a diffeomorphism of $T^8$ defined by a matrix $A\in SL(8, \mathbb{R})$. We choose $A$ as in the following Lemma:

\hfill

 \begin{lemma}\label{pseudo-hyperkahler}
Consider on $\R^8$  the hypercomplex structure $(I, J, K)$ defined, in terms of the standard basis $(e_1, \ldots, e_8)$,  by
$$
\begin{array}{llll}
I e_1 = e_3, & I e_2 = e_4,&  I e_5 = - e_7, & I e_6 = - e_8,\\
J e_1 = - e_5, &  J e_2 = - e_6,  &  J e_3 = - e_7 &  J e_4 = - e_8,
\end{array}
$$
and the pseudo hyperhermitian metric
$$
h= (e^1)^2 - (e^2)^2 +(e^3)^2 - (e^4)^2+ (e^5)^2 - (e^6)^2 +   (e^7)^2 - (e^8)^2.
$$
 Then there exists an integer matrix $$A =  \left( \begin{array}{ccccccccc}  1 & 0 & 1 &  0 &  -1 &  -1 &  0 &   1\\[2pt]   0 &  -1 &  0 &  -1 &  -1 &  0 &  1 &  1\\[2pt]  -1 &  0&   1 &  0 &  0 &  1 &  1 &  1\\[2pt]   0 &  1 &  0 &  -1 &  1 &  1 &  1 &  0\\[2pt]    1 &  1 &  0 &  -1 &  1 &  0 &  1 &  0\\[2pt]   1 &  0 &  -1 &  -1 &  0 &  -1 &  0 &  -1\\[2pt]   0 &  -1 &  -1 &  -1 &  -1 &  0 &  1 &  0\\[2pt]   -1 &  -1 &  -1 &  0 &  0 &  1 &  0 &  -1 \end{array}  \right ) \in SL(8, \Z)$$  preserving the   pseudo hyperhermitian structure $(I, J,K, h)$ and such  $A^k \neq Id$,  for every $k >  1$.
\end{lemma}

\hfill

{\bf Proof:}   Since $A$ commutes with the matrices associated to $I$ and $J$ with respect to the standard basis $(e_1, \ldots, e_8)$ and $A^t H A = H$, where $H$ is the matrix associated to $h$, we have that
$A$ preserves the  pseudo hyperhermitian structure $(I, J,K, h)$. Moreover,   $A^k \neq Id$, for every $k >  1$,  since $A$ has eigenvalues
$\pm i (1+\sqrt{2}) $ and  $ \pm  i(\sqrt{2}-1)$ of multiplicity two.
\endproof

\hfill

Recall that an indefinite or pseudo hyperhermitian metric is called {\bf pseudo-hyperk\"ahler} if the corresponding fundamental forms are closed.
In particular, the Lemma claims that $A$ preserves the pseudo-hyperk\"ahler structure on $\mathbb{R}^8$  and the lattice generated by $e_1,..,e_8$. Such a matrix defines a hyperbolic diffeomorphism $f_A$ of $T^8$ which preserves the pseudo-hyperkahler structure and in particular both the hypercomplex structure and the holomorphic symplectic form. We note that any such matrix which preserves a positive-definite metric, has all eigenvalues on the unit circle so cannot be hyperbolic.  To formulate the next Theorem, recall that for every hypercomplex manifold there exists a unique torsion free connection preserving the hypercomplex structure, which is called the {\bf Obata connection}.

\hfill

\begin{theorem} \label{pseudohkex}
Let $f_A$ be  the diffeomorphism of $T^8$ defined by the matrix $A$ and  $S = S(T^8, f_A, id, id, id)$ the hyperbolic toric suspension, where $T^8$ denotes the 8-dimensional real torus obtained as a quotient of  $\mathbb{R}^8$ by the standard lattice.
 Then  $S$ admits a pseudo-hyperk\"ahler structure $(I, J, K, h).$ In particular for the  complex structure $I$,  $S$  carries both hypercomplex and holomorphic symplectic structures, but no K\"ahler metrics. Moreover, the Obata connection of the  hypercomplex structure is flat.
\end{theorem}

{\bf Proof:}   Since $f_A$ preserves both the standard  hypercomplex and  holomorphic symplectic structures on $\mathbb{R}^8$,  the previous  structures descend to  $T^8$. Then the action $(p, t)   \rightarrow  (f_A^{-n}(p), t + n)$ of $\mathbb{Z}$ on $T^8\times \mathbb{R}^4$ preserves the induced natural structures obtained as  product of the ones on $T^8$ and the canonical hyperk\"ahler structure on $\mathbb{R}^4$.  Moreover,  the   hypercomplex structure on $T^8\times \mathbb{R}^4$  is clearly compatible  with  flat Obata connection. As a consequence  all structures descend to the quotient $S$.  The fact that $S$ is non-K\"ahler follows  from  \ref{_suspension_non-Kahler_Theorem_}, because  $f_A$ does not preserve any K\"ahler class.\endproof

\hfill

{\remark Note that pseudo-hyperk\"ahler structures on  $12$-dimensional compact solvmanifolds are constructed in \cite{Yamada} and are associated to the almost abelian  Lie algebras $$\Psi_I (\Psi_J (\frak g)) = {\rm{span}}_{\R} \{ U_1^1, U_1^2, U_1^3, U_1^4, V_1^1, V_2^1,  V_1^2, V_2^2, V_1^3, V_2^3, V_1^4, V_2^4 \}$$ with Lie bracket
$$
[U_1^1, V_j^h] = c_{1j}^1 V_1^h +  c_{1j}^2 V_2^h, \quad j =1,2,  \, \,  h=1,2,3,4,
$$
and hypercomplex structure $(I, J, K)$ defined by
$$
\begin{array}{llllll}
I U_1^1 = U_1^2, & I V_1^1 = V_1^2,  & IV_2^1 = V_2^2, &  I U_1^4 = U_1^3, & I V_1^4 = V_1^3, & I V_2^4 = V_2^3,\\[3pt]
J U_1^1= U_1^3, &  J V_1^1 = V_1^3,  & J V_2^1 = V_2^3, &  J U_1^2 = U_1^4, & I V_1^2 = V_1^4, &  I V_2^2 = V_2^4.
\end{array}
$$
So in particular the associated solvable Lie group is a  semidirect product of the form $(\R \ltimes \R^8) \times \R^3$.
In the notation of  \cite{Yamada} $\frak g = \frak a \ltimes \frak b$, with $\frak a = {\rm{span}}_{\R} \{ U_1^1 \}$ and  $\frak b =  {\rm{span}}_{\R}\{ V_1 ^1, V_2^1 \}$ and by Theorem 6.6  in \cite{Yamada}  if $\mathfrak{b}$ has a non-degenerate $2$-form which is closed on $\frak g$, then $($$\Psi_I (\Psi_J (\frak g)), I, J,K)$ admits a compatible pseudo-hyperk\"ahler structure. The   previous condition  is satisfied if $c_{11}^1 = -c_{12}^2$.  The hyperbolic toric suspension $S = S(T^8, f_A, id, id, id)$ corresponds to the  compact solvmanifold constructed as a quotient of the solvable  Lie group  $H$  whose Lie algebra is $ \frak h := \Psi_I (\Psi_J (\frak g))$ with  $ad_{U_1^1}= {\text {diag}} (1,-1, 1,-1, 1,-1,1,-1)$. If we  consider the basis
$$
\begin{array}{llllll}
f_1 = V_1^1, & f_2 = V_2^1 & f_3 = V_1^2, & f_4 = V_2^2, &  f_5 = V_1^3, & f_6 = V_2^3,\\
 f_7 = V_1^4, &  f_8 = V_2^4, &  f_9 = U_1^1, &  f_{10} = U_1^2, &  f_{11} = U_1^3, &  f_{12} = U_1^4,
 \end{array}
$$
the structure equations  of $\frak h$ are
\begin{equation}\label{structureq}
\begin{array}{l}
d f^i= f^i\wedge f^9, \,  i= 1,3,5,7, \quad d f^j = - f^j \wedge f^9, \,  j =2,4,6,7, \\[2pt]
 d f^k =0, k = 9,\ldots 12.
 \end{array}
\end{equation}
The   pseudo-hyperk\"ahler    structure on $\frak h$  is given by  $(I, J, K,  \omega_I, \omega_J, \omega_K)$, where
 $$
 \begin{array}{l}
 \omega_I =  2 (-f^1 \wedge f^2 - f^3 \wedge f^4 + f^5 \wedge  f^6  +f^7 \wedge f^8 + f^9 \wedge f^{10} - f^{11} \wedge f^{12}),\\[2pt]
 \omega_J =  2 (f^1 \wedge f^8 + f^4 \wedge f^5 - f^2 \wedge  f^7  - f^3 \wedge f^6 + f^9 \wedge f^{11} + f^{10} \wedge f^{12}),\\[2pt]
  \omega_K=  2 (f^1 \wedge f^6 - f^4 \wedge f^7 - f^2 \wedge  f^5  +f^3 \wedge f^8 - f^9 \wedge f^{12} + f^{10} \wedge f^{11}).
 \end{array}
 $$
 With respect to the basis of $(1,0)$-forms with respect to $I$
\begin{equation}\label{etai}
 \begin{array}{l}
 \eta_1 = f^1 + i f^3, \, \eta_2 = f^2 + i  f^4, \,  \eta_3 = f^5 - i f^7, \\[2pt]
  \eta_4 = f^6 - i f^8, \, \eta_5 = f^9 + i f^{10}, \, \eta_6 = f^{11} - f^{12}
  \end{array}
\end{equation}
 we have
\begin{equation}\label{Jetai}
 J \eta_1 = \overline \eta_3, \, J \eta_2 = \overline \eta_4, \, J \eta_3 = -  \overline \eta_1, \, J \eta_4 = - \overline \eta_2, \, J \eta_5 = \overline \eta_6, \, J \eta_6 = - \overline \eta_5
\end{equation}
 and the associated $(2,0)$-form $\omega_J + i \omega_K $ is given by
 $$ \omega_J + i \omega_K= 2  (\eta_5 \wedge \eta_6 + i \eta_1 \wedge \eta_4 - i \eta_2 \wedge \eta_3).
 $$
Note that
 $$
  J(\omega_J + i \omega_K ) =  2  ( \overline \eta_5 \wedge  \overline \eta_6 - i  \overline \eta_1 \wedge  \overline \eta_4 + i  \overline \eta_2 \wedge  \overline \eta_3)
 $$
and that the  two   $(4,0)$-forms $ \eta_1 \wedge \eta_3 \wedge \eta_5 \wedge \eta_6$ and $ \eta_2 \wedge \eta_4 \wedge \eta_5 \wedge \eta_6$ are both $\partial$-exact.
 }
\hfill


\subsection{Compatible metric structures}
\label{_HKT_Section_}

On a hypercomplex manifold $(M, I, J, K)$ we always have a {\bf hyper-Hermitian} (positive-definite) metric, that is  a metric compatible with the complex structures $I, J, K$. When the fundamental forms $\omega_I, \omega_J, \omega_K$ are closed the metric is hyperk\"ahler, but in the example  constructed in  \ref{pseudohkex}  such metric doesn't exist.  A  generalization of hyperk\"ahler condition is the condition $\partial \Omega =  0$, where $\Omega = \omega_J+i\omega_K$,  in which case the metric is called  {\bf hyperk\"ahler  with torsion} (shortly HKT) \cite{GP}. In fact one can characterize the HKT condition in terms of $\Omega$: if there is a $(2,0)$ form  $\Omega$ with respect to  $I$, such that $\partial\Omega =0$, $\Omega(JX,JY) = -\overline{\Omega(JX,JY),}$ and $\Omega(X, J\overline{X}) >0$ for  every non-zero $(1,0)$  vector field $X$, then the metric  defined by $g(X,Y) = Re\Omega(X,J\overline{Y})$ is HKT.  The HKT metric is a good candidate for a quaternionic analog of K\"ahler metrics in complex geometry - it arises from a local quaternionic-subharmonic potential and gives rise to a Hodge theory (see \cite{GP} and \cite{V}).  The existence or non-existence of HKT metrics, in 8-dimensional case depends on purely holomorphic data  (see \cite{GLV}). In hypercomplex geometry the analog of the balanced condition is called {\bf quaternionic balanced} (see \cite{LW}) and such metric satisfies $\partial(\Omega^{n-1})=0$, where $2n$ is the complex dimension of the manifold.  We have the following:

\hfill

\begin{theorem}
The hypercomplex manifold  $S = S(T^8, f_A, id, id, id)$ from \ref{pseudohkex}  admits a quaternionic balanced metric, but admits   no HKT metrics.
\end{theorem}

\hfill

{\bf Proof:} We consider the solvmanifold model  of $S$ from Remark 7.3.  We use the same $(1,0)$-forms $\eta_i, 1\leq i \leq 6$,  and  complex structure $J$  as in \eqref{etai} and \eqref{Jetai}. Note that the hypercomplex structure has Obata holonomy in $SL(n,  \mathbb{H})$, since it has a  closed and real $(6,0)$-form.  From the structure equations \eqref{structureq} we see that the $(2,0)$-form form $\Omega = \eta_1\wedge \eta_3 + \eta_2\wedge \eta_4+ \eta_5\wedge \eta_6$ satisfies the condition $$\partial\Omega^2 = \partial (\eta_1 \wedge \eta_3 \wedge \eta_5 \wedge \eta_6 + \eta_2 \wedge \eta_4 \wedge \eta_5 \wedge \eta_6 + \eta_1 \wedge \eta_2 \wedge \eta_3 \wedge \eta_4) = 0,$$ but $\partial\Omega \neq 0$, so $\Omega$ defines a  quaternionic-balanced metric. On the other side, by averaging argument (see \cite{FG}), if there is any HKT metric, then there exists an invariant one. Working by contradiction, we assume that there is  a $(2,0)$-form $\tilde \Omega$ with $\partial(\tilde{\Omega}) = 0$ which is $J$-anti-invariant and positive, so defines an HKT metric. Then $\tilde{\Omega}$ has the form $$\tilde{\Omega} = \sum_{\alpha, \beta} a_{\alpha\overline{\beta}}  \, \eta_\alpha\wedge J\overline{\eta_\beta},$$ where  $a_{\alpha\overline{\beta}}$ is a Hermitian and positive definite matrix. Now we can adapt the Harvey-Lawson property for $SL(n, \mathbb{H})$ manifolds from \cite{GLV} and use it explicitly. Since the $(4,0)$-form $\alpha= \eta_1 \wedge \eta_3 \wedge \eta_5 \wedge \eta_6 + \eta_2 \wedge \eta_4 \wedge \eta_5 \wedge \eta_6$ is $\partial$-exact and the $(6,0)$-form $\beta = \eta_1\wedge\eta_2\wedge\eta_3\wedge\eta_4\wedge\eta_5\wedge\eta_6$ is closed, then we have $$\int_M \tilde{\Omega}\wedge\alpha\wedge \overline{\beta} = 0$$ by integration by parts. On the other side $vol = \beta\wedge\overline{\beta}$ is a volume form  and $\tilde{\Omega}\wedge\alpha\wedge \overline{\beta} = (a_{1 \overline 1} + a_{2 \overline 2}) vol   >0$, so $$ \int_M \tilde{\Omega}\wedge \alpha \wedge \overline{\beta} > 0$$
and we get a contradiction. \endproof

\section{Pluriclosed metrics from suspensions} \label{Sect7}

We recall the following
\hfill

\definition
Let $(M,I)$ be a  complex  manifold.
We say that  an $I$-Hermitian metric $h$ is {\bf pluriclosed}  (or {\bf SKT}) if  its fundamental form is  $\partial \overline \partial$-closed.
\hfill

Using automorphisms of K\"ahler manifolds it is  also possible to construct suspensions admitting pluriclosed metrics in the following way.

Let $(X^3, \xi, \eta,  \varphi, \Phi) $ be  a 3-dimensional  Sasakian  manifold
and   $(Y^{2n}, I, g,  \omega)$  be a  K\"ahler manifold of complex dimension $n$.

 On  the product $X^3 \times Y^{2n}  \times \R$    we can define a complex structure $\tilde I$  such that
$\tilde I (\xi) =\frac{\partial}{ \partial_t}$, where $t$ is the coordinate on $\R$ and $ \tilde I= I$  on $Y^{2n}$.
The $2$-form
$$
\tilde \omega = \eta \wedge dt + d \eta  +  \omega
$$
 is then a   positive  $(1,1)$-form on $(X^3 \times Y^{2n}  \times \R, \tilde I)$.
Since   $d \eta = \Phi$, by a direct computation we obtain
 $$
d \tilde \omega  = \Phi \wedge dt
$$
and
$$
d  T^B =  d (\tilde I  d \tilde \omega) = -  d( \Phi \wedge \eta) =0,
$$
where $T^B = \tilde I d  \tilde  \omega$ is the so-called {\bf Bismut torsion form}. Note that $ \Phi \wedge \eta$ is closed since it is a $3$-form on $X^3$.
Therefore we have the following

\hfill

\begin{theorem}  Let $(X^3, \xi, \eta,  \varphi, \Phi)$ be  a Sasakian 3-dimensional manifold,
 $(Y^{2n}, I, g,  \omega)$  a  K\"ahler manifold and $f = (f_1,f_2)$ a diffeomorphism   of $X^3 \times Y^{2n}$ such that
 $f_1$ is a diffeomorphism of $X^3$  preserving the Sasakian structure $(\xi , \eta,  \varphi, \Phi)$ and $f_2$ is a diffeomorphism of $Y^{2n}$  preserving the K\"ahler structure $(I,  g, \omega)$. Then the suspension of  $X^3 \times Y^{2n}$  by $f$ is non-K\"ahler and has a pluriclosed   metric.
\end{theorem}

\hfill

To prove  that the suspension of   $(X^3, \xi, \eta,  \varphi, \Phi)$ by $f$  is non-K\"ahler  we can use the   Harvey-Lawson characterization of non-K\"ahler manifolds \cite{HL} and the observation that $d\eta=\Phi$ is a  non-zero (weakly) positive and exact  $(1,1)$-form.

\hfill

\begin{example} An application of the previous construction  gives the example  of compact solvmanifold constructed in \cite{_Fino:Paradiso}.    More precisely, let $G$ be the simply connected $3$-step  solvable Lie group with structure equations$$
\left \{ \begin{array}{l}
de^1 = e^2 \wedge e^3,\\
de^2 = -e^2 \wedge e^8,\\
de^3 = e^3 \wedge e^8,\\
de^4 = b \, e^5 \wedge e^8,\\
de^5 = -b \, e^4 \wedge e^8,\\
de^6 = b \, e^7 \wedge e^8\\
de^7 = -b \, e^6 \wedge e^8,\\
de^8 = 0,
\end{array}
\right.
$$
with $b= \frac{2 \pi}{\log(2+\sqrt{3})}$. By \cite{_Fino:Paradiso} $G$ has the left-invariant complex structure
$$
\begin{array}{l}
I e_1 =-  e_2,  \,  I e_3 =e_8,  \, I e_4 = e_5, \,  I  e_6 = e_7,
\end{array}
$$
 and  admits a compact quotient by a lattice $\Gamma$.
The $I$-Hermitian metric  $g = \sum_{i =1}^{8} (e^i)^2$  is pluriclosed since
the Bismut torsion $3$-form $T^B = I d \omega$ is  the closed $3$-form $- e^1 \wedge e^2 \wedge e^3$. The compact solvmanifold  $\Gamma \backslash G$  can be  obtained   as  a suspension of the product of the 3-Sasakian manifold given by the compact quotient of the real $3$-dimensional Heisenberg group by a lattice and the standard   torus  $T^4$. Moreover, the compact solvmanifold   can be viewed also as  the total space of a bundle over a circle with fibre   a circle bundle over a 6-torus.

\end{example}

\hfill

{\bf Acknowledgements:} We thank the anonymous referee for the careful reading of the manuscript and the numerous suggestions and clarifications which led to a considerably improved exposition.

\hfill

{\bf Conflict of interest statement.}  On behalf of all authors, the corresponding author states that there is no conflict of interest.

\medskip

\small

\noindent
{\sc Anna Fino}\\
{\sc Dipartimento di Matematica G. Peano\\
Universit\'a di Torino\\
via  Carlo Alberto 10, 10123 Torino, Italy}\\
Also:\\
{\sc  Department of Mathematics and Statistics\\
 Florida International University\\
Miami Florida, 33199, USA\\
\tt  annamaria.fino@unito.it, afino@fiu.edu}\\

\noindent {\sc Gueo Grantcharov\\
{\sc Department of Mathematics and Statistics\\
Florida International University\\
Miami Florida, 33199, USA}\\
\tt grantchg@fiu.edu}\\

\noindent {\sc Misha Verbitsky}\\
    {\sc Instituto Nacional de Matem\'atica Pura e
            Aplicada (IMPA) \\ Estrada Dona Castorina, 110\\
Jardim Bot\^anico, CEP 22460-320\\
Rio de Janeiro, RJ - Brasil}\\
{\tt verbit@impa.br}\\


\small

\begin{thebibliography}{AV1}

\bibitem[AFLM]{AFLM} de Andr\'es, L.C.,  Fern\'andez, M.,  de Le\'on, M.  Menc\'a, J.J.,
{\em Some six dimensional compact symplectic and complex
  solvmanifolds}, Rendiconti di Mat. Roma {\bf 12} (1992),
59--67.


\bibitem[AB]{AB} Alessandrini, L., Bassanelli, G.,
{\em Transforms of currents by modifications and
  $1$-convex manifolds}, Osaka J. Math. {\bf 40} (2003),
717--740.

\bibitem[AB2]{AB2} Alessandrini, L., Bassanelli, G.,
{\em Wedge product of positive product of positive
  currents and balanced manifolds}, Tohoku Math. J. {\bf
  60} (2008), 123-134.

\bibitem[AV]{_AV:Aut_}
Amerik, E., Verbitsky, M.,
{\em Construction of automorphisms of hyperk\"ahler manifolds,}
Compos. Math.  {\bf 153} (2017), 1610--1621.

\bibitem[Bea]{_Beauville_}
Beauville, A., {\em Varietes K\"ahleriennes dont la premi\`ere classe de
Chern est nulle.}  J. Diff. Geom. {\bf 18} (1983),  755--782.

\bibitem[Be]{Beri}  Beri, P.,   On birational
  transformations and automorphisms of some hyperk\"ahler
  manifolds, PhD-thesis,  Universit\'e de Poitiers, 2020.



\bibitem[Bes]{_Besse:Einst_Manifo_}
Besse,
A., {\em Einstein Manifolds}, Springer-Verlag, New York (1987).

\bibitem[Bo]{_Bogomolov_}
Bogomolov, F. A., {\em On the
  decomposition of K\"ahler manifolds with trivial
canonical class,} Math. USSR-Sb.  {\bf 22}  (1974),  580--583.

\bibitem[BKLV]{_BKLV_} Bogomolov, F., Kamenova, L.,  Lu, S.,
Verbitsky, M., {\em On the Kobayashi pseudometric, complex
automorphisms and hyperk\"ahler manifolds}, Geometry over
nonclosed fields, 1-17, Simons Symp., Springer, Cham,
2017.

\bibitem[Ca]{_Cantat:K3_}
Cantat, S.,
{\em Dynamique des automorphismes des surfaces $K3$},
Acta Math. {\bf 187}  (2001), no. 1, 1--57.

\bibitem[Ch]{C14} Chiose,  I., {\em Obstructions to the existence of K\"ahler structures on compact complex manifolds}, Proc. Amer. Math. Soc.  {\bf 142}, No. 10 (2014), 3561--3568.

\bibitem[Dem]{_Demailly:agbook_}
  Demailly, J.-P.,  Complex analytic and differential
    geometry, 2012   (\url{http://www-fourier.ujf-grenoble.fr/*demailly/manuscripts/agbook.pdf}).

\bibitem [DS]{DinhSibony10}
Dinh, T.-C., Sibony, N.,  {\em Super-potentials for currents on compact K\"ahler manifolds and dynamics of automorphisms},
J. Algebraic Geom. {\bf  19} (2010), no. 3, 473--529.


\bibitem[Fe]{Fe}  Fei, T., {\em Construction of Non-K\"ahler Calabi-Yau Manifolds and new solutions to the Strominger
System}, Adv. Math. {\bf 302}  (2016) 529--550.

\bibitem[FMS]{FMS}  Fern\'andez, M.,   Mun\v oz, V., Santisteban, J.  A.,
{\em Cohomologically K\"ahler manifolds with no K\"ahler metrics},
Int. J. Math. Math. Sci. {\bf 52} (2003),  3315--3325.

\bibitem[FG] {FG}  Fino, A.,  Grantcharov, G., {\em Properties of manifolds with skew-symmetric torsion and special holonomy}, Adv. Math. {\bf 189} (2004), no. 2, 439--450.

\bibitem[FGV]{FGV} Fino, A., Grantcharov, G., Vezzoni, L., {\em Astheno-K\"ahler and balanced structures on fibrations}, Int. Math. Res. Not. IMRN {\bf 2019}, no. 22, 7093--7117.


\bibitem[FP]{FP}  Fino, A., Paradiso, F., {\em  Balanced Hermitian structures on almost abelian Lie algebras}, J. Pure Applied Algebra {\bf 227} (2023),  no. 2, Paper No. 107186.


\bibitem[FP2]{_Fino:Paradiso}  Fino, A., Paradiso, F., {\em  Hermitian structures on a class of almost nilpotent solvmanifolds}, J.  Algebra
{\bf 609}  (2022), 861--925.

\bibitem[FS]{_Freibert-Swann_}  Freibert, M., Swann, A.,
  {\em Compatibility of balanced and SKT metrics on
    two-step solvable Lie groups}, arXiv:2203.16638.

\bibitem[Fu]{_Fujiki_} Fujiki, A., {\em On the de Rham Cohomology Group of a Compact K\"ahler Symplectic
Manifold}, Adv. Stud. Pure Math. {\bf 10} (1987),  105--165.

\bibitem[FV]{FV16} Fino,  A.,Vezzoni,  L., {\em On the existence of balanced and SKT metrics on nilmanifolds}, Proc. Amer. Math. Soc., {\bf 144}(6), (2016), 2455--2459.

\bibitem[FLY]{FLY12}   Fu, J.,  Li ,  J., Yau, S.-T., {\em Balanced metrics on non-K\"ahler Calabi-Yau threefolds}, J. Diff. Geom. {\bf 90} (2012), 81--129.

\bibitem[GiPo]{GiPo} Giusti, F., Podest\'a,  F., {\em Real semisimple Lie groups and balanced metrics},  Rev. Mat. Iberoam. {\bf 39} (2023), no.2, 711--729.
    
\bibitem[Gall]{Gallier:diffgeom_and_Lie} Gallier, J., Quaintance, J., {\em  Differential Geometry and Lie Groups: A computational approach},  Springer (2020).


\bibitem[Gol]{_Goldman_}  Goldman,  W.,
Geometric structures on manifolds, 2021   (\url{http://www.math.umd.edu/*wmg/gstom.pdf  }).

\bibitem[GLV]{GLV}  Grantcharov, G., Lejmi, M., Verbitsky,
  M., {\em  Existence of HKT metrics on hypercomplex
    manifolds of real dimension 8},  Adv. Math. {\bf 320}
  (2017), 1135--1157.

\bibitem[GP]{GP} Grantcharov, G., Poon, Y. S.,  {\em Geometry of hyperK\"ahler connections with torsion}, Commun. Math. Phys. {\bf 213}  (2000), 19--37.

\bibitem[Gu1]{Gua1}  Gualtieri,  M.,  Generalized complex geometry, PhD-thesis, University of Oxford, 2003, arXiv:math/0401221.

\bibitem[Gu2]{Gua2}   Gualtieri, M.,  {\em Generalized  K\"ahler geometry}, Comm. Math. Phys. {\bf 331} (2014), no. 1, 297--331.

\bibitem[HL]{HL} Harvey, R.,  Lawson Jr.,  B. H.,   {\em An intrinsic characterization of K\"ahler manifolds}, Invent. Math. {\bf 74} (1983), no. 2, 169--198.

\bibitem[Hi]{Hi}  Hitchin, H., {\em Generalized Calabi-Yau Manifolds}, Q. J. Math. {\bf 54} (2003), no. 3, 281--308.

\bibitem[Hu]{Hu} Hull, C.,  {\em Superstring compactifications with torsion and space-time supersymmetry}, In Turin 1985 Proceedings “Superunification and Extra Dimensions” (1986), 347–375.


\bibitem[Ka]{_Kapovich:Kleinian_}
Kapovich,  M.,  Kleinian groups in higher dimensions. In
"Geometry and Dynamics of Groups and Spaces. In memory of
Alexander Reznikov'', M.Kapranov et al (eds). Birkhauser,
Progress in Mathematics, Vol. 265, 2007, p. 485-562.

\bibitem[LU]{LU} Latorre, A., Ugarte, L.,  {\em On non-K\"ahler compact complex manifolds with balanced and astheno-K\"ahler metrics},  C. R. Math. Acad. Sci. Paris {\bf 355}, no. 1 (2017),  90--93.

\bibitem[LW]{LW}  Lejmi, M.,  Weber, P., {\em Quaternionic Bott-Chern cohomology and existence of HKT metrics}, Q. J. Math. 68 (2017), no. 3, 705--728.

\bibitem[Li]{Li}  Li, H., {\em Topology of co-symplectic/co-K\"ahler manifolds},  Asian  J. Math. {\bf 12}  (2008),   527--544.

\bibitem[Ma]{_Manjarin_} Manjar\'in, M.,  {\em Normal almost contact structures and non-K\"ahler compact complex manifolds},  Indiana Univ. Math. J. {\bf 57} (2008), no. 1, 481--507.

\bibitem[Mic]{_Michelson_}
Michelsohn, M. L., {\em On the existence of special metrics in complex geometry}, Acta Math. {\bf 143}  (1983), 261--295.

\bibitem[OV]{OV} Ornea, L., Verbitsky, M., {\em Principles of Locally Conformally K\"ahler Geometry},  arXiv:2208.07188.

\bibitem[O]{Otiman} Otiman, A., {\em Special Hermitian metrics on  Oeljeklaus-Toma manifolds}, Bull. Lond. Math.  Soc. {\bf 54} (2022), 655--667.


\bibitem[St]{St}   Strominger, A. E.,  {\em Superstrings with torsion}, Nuclear Phys. B {\bf 274}(2) (1986), 253--284.


\bibitem[QW]{QW}  Qin, L:,  Wang, B.,
{\em A family of compact complex and symplectic Calabi-Yau manifolds that are non-K\"ahler},
Geom. Topol. {\bf 22}  (2018), 2115--2144.


\bibitem[STW]{STW} Sz\'ekelyhidi, G., Tosatti, V.,   Weinkove, B., {\em Gauduchon metrics with prescribed volume form},  Acta Math. {\bf 219} (2017), no. 1, 181--211.


\bibitem[Ve02] {V} Verbitsky,  M.,  {\it HyperKähler manifolds with torsion, supersymmetry and Hodge theory}, Asian J. Math. {\bf 6} (2002), no. 4, 679–712.


\bibitem[Ve14]{Ve14} Verbitsky,  M., {\it Rational curves and special metrics on twistor spaces}, Geom. Topol. {\bf 18}  (2014), no. 2, 897--909.

\bibitem[Ya]{Yamada} Yamada, T.,
{\em A construction of compact pseudo-K\"ahler
  solvmanifolds with no K\"ahler structures}, Tsukuba
J. Math. {\bf 29} (2005),  79-109.

\end{thebibliography}
\end{document}